\documentclass{llncs}

\usepackage{llncsdoc}
\usepackage{graphicx}
\usepackage{fourier}

\begin {document}
\title{Two-Fold Circle-Covering of the Plane under
Congruent Voronoi Polygon Conditions}
\titlerunning{Two-Coverage}

\author{Jingchao Chen}
\institute{School of Informatics, Donghua University \\
2999 North Renmin Road, Songjiang District, Shanghai 201620, P. R.
China \email{chen-jc@dhu.edu.cn}}

\maketitle

\begin{abstract}
The $k$-coverage problem is to find the minimum number of disks
such that each point in a given plane is covered by at least $k$
disks. Under unit disk condition, when $k$=1, this problem has
been solved by Kershner in 1939. However, when $k > 1$, it becomes
extremely difficult. One tried to tackle this problem with
different restrictions. In this paper, we restrict ourself to
congruent Voronoi polygon, and prove the minimum density of the
two-coverage with such a restriction. Our proof is simpler and more
rigorous than that given recently by Yun et al.~\cite{Yun:2013,Yun:2011}.
\end{abstract}

\section{Introduction}

A set of circles is said to form a $k$-fold covering if every point
of the plane belongs to at least $k$ circles. The problem of finding
the thinnest (or minimum density) $k$-fold covering has been studied
for a long time.

In 1939, Kershner~\cite{Kershner:1939} solved the one-coverage
problem with equal disks, and proved that the triangle lattice
pattern is the best, and its minimum density $\vartheta$ is
$\vartheta$=$\pi / A_6$, where $A_6 = 3\sqrt{3}/2$ is the area of a
regular hexagon inscribed in a unit circle. Here the density of covering may be
defined simply as the ratio of the total area of disks covering a given region $M$
to the area of $M$.

The $k$-coverage problem
with $k > 1$ is much more difficult than the one-coverage problem. So far
only a few special $k$-coverage problems were solved.

For example, in 1957, Blundon~\cite{Blundon:1957} restricted himself to lattice-coverings, i.e., suppose that
the centers of the circles form a lattice. Let $\vartheta^{(k)}$ be
the minimum density of $k$-fold lattice covering by equal circles.
Blundon proved that $\vartheta^{(2)} = 2 \vartheta$,
$\vartheta^{(3)} = 2.841 \vartheta$, $\vartheta^{(4)} = 3.608
\vartheta$, $\ldots$.

In 1960, instead of Blundon's single lattice covering, Danzer~\cite{Danzer:1960} considered a  multiple lattice-covering, and
obtained the minimum density $D_2$ of the two-coverage: $2.094\ldots
\leq D_2 \leq 2.347\ldots$.

In 1976, Toth~\cite{Toth:1976} used the notion of the $k$-th
Dirichlet (Voronoi) cell to estimate the the lower bound of the
$k$-coverage. Let $D_k$ be the minimum density of the $k$-coverage.
Toth proved that  $D_k \geq \frac{\pi}{3}\csc \frac{\pi}{3k}$, where $\csc \theta = \frac{1}{\sin \theta}$.

Recently Yun et al. \cite{Yun:2013,Yun:2011} investigated congruent
Voronoi polygon coverings, i.e., suppose that Voronoi polygons
generated by the centres of the equal circles are congruent. Then
under such a restriction, how large is the minimum density of the
two-coverage? To answer this problem, we formulate the notion of
Voronoi polygon as follows.

\begin{definition} \texttt {(Voronoi polygon)}
Let $a_1, a_2, \ldots, a_n$ be $n$ points in a Euclidean plane $P$.
 The Voronoi polygon $V(a_i)$ ($1 \leq i \leq n$) is defined as the set
of all points in $P$, which are closer to the center $a_i$ than to
any other center $a_j$, i.e., $V(a_i)=\{ x | d(x,a_i) \leq d(x,a_j)$, $1 \leq j \leq n$\},
where $d(\bullet,\bullet)$ denotes the Euclidean distance.
\end{definition}

Using the notion of Voronoi polygon, the covering density may defined as follows.

\begin{definition} \texttt {(covering density)}
The covering density is the ratio of the total area of disks to
the total area of Voronoi polygons generated by the centers of the disks.
\end{definition}

Under congruent Voronoi polygon conditions, the
covering density is actually the ratio of the area of a disk to
the area of Voronoi polygon inscribed in the disk.
Yun et al.~\cite{Yun:2013,Yun:2011} claimed that the minimum density of the
two-coverage under congruent Voronoi polygon conditions is
2$\vartheta$. However, their proof on it is considered not to be
rigorous. For example, Lemma 4.1 in \cite{Yun:2013} is not
straightforward. However, Yun et al. did not prove it. That is, they
did not show why there are only 12 types of congruent Voronoi
polygonal tessellations. In \cite {Branko:1978,Branko:1987},
Gr\"{u}nbaum et al. listed the 107 polygonal isohedral types of
tilings. Yun et al. did not show how to reduce the 107 types to 12
types and exclude non-isohedral tilings (Voronoi polygon is not necessarily
isohedral). Also, in the proof of Lemma 4.7 of \cite{Yun:2013},
the authors proved type ($h$), using the unproven inequality:
$2\sqrt{r^2-(\frac{a^2}{\sqrt{a^2+1}})^2} \geq\sqrt{a^2+1}$. We
cannot conclude that this inequality holds.

  This paper addresses the same two-coverage problem as
  Ref.~\cite{Yun:2013}, i.e., two-fold covering with congruent Voronoi polygonal restrictions.
  However, we adopt a different proving strategy to provide a simpler and more rigorous
  proof than Yun et al.

\section{Minimum Density of the Two-Coverage
with Congruent Voronoi Polygon Restrictions}
\begin{figure}
\centering
\includegraphics[height=4.0cm]{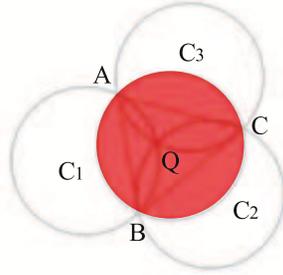}
\caption{Three intersections $A, B, C$ of three equal circles $C_1,
C_2, C_3$ lie on the red circle of the same radius.}
\label{triangleA}
\end{figure}


The main theorem of this paper will be presented after the following
two lemmas.

\begin{lemma} $\mathrm {(Johnson's~~theorem)}$
As shown in Figure \ref{triangleA}, suppose that three equal circles
$C_1, C_2, C_3$ with radius $r$ intersect at a point $Q$, and three other intersections
are $A$, $B$ and $C$. Then the radius of
circumcircle of triangle $ABC$ is $r$ also, i.e., these four circles
are equal.
\end{lemma}
%

\begin{proof}
This lemma is called Johnson's theorem proved by Johnson in 1916 \cite{Johnson:1929,Johnson:1916}. We prove briefly it as follows.
  
$\angle AQB$ subtends arc $\widearc {AB}$ of circle $C_1$. If we express them in radians,
we have $\angle AQB = \frac{\widearc {AB}}{2} = \frac{2\pi - \widearc {AQB}}{2}$. Similarly, we have have $\angle BQC
=  \frac {2\pi - \widearc {BQC}}{2}$ and $\angle CQA
=  \frac {2\pi - \widearc {CQB}}{2}$. Clearly $\angle AQB + \angle BQC + \angle CQA = 2\pi$. Therefore, we have
$\widearc {AQB} + \widearc {BQC} + \widearc {CQA}  = 2\pi$.
In addition, clearly, for the red circle, we have $\widearc {AB} + \widearc {BC} + \widearc {CA} = 2\pi$.
Suppose that the radius of the red circle is not equal to the radius of $C_1, C_2, C_3$. Then it
implies $\widearc {AB} + \widearc {BC} + \widearc {CA} \neq 2\pi$, since the same
chord on the circles of different radii subtends the arc with
different angles. This is a contradiction.
\end{proof}

\begin{lemma} 
Under congruent Voronoi triangle conditions, the minimum density of
the two-coverage with equal circles is 2$\vartheta$, where
$\vartheta$ is the minimum density of the one-coverage.
\end{lemma}

\begin{proof}
To satisfy the conditions of the lemma, i.e., Voronoi triangles
inscribed in each circle are equal and achieve two-coverage, the
tessellation pattern must be the one shown in Figure
\ref{triangleA}. In this case, the density is the ratio of the area
of a circle to the area of $\triangle ABC$. It is well known that
the area of an inscribed triangle reaches its maximum when it is
regular. Therefore, the minimum density is $\frac{4\pi}{\sqrt{27}}$
= 2$\vartheta$.
\end{proof}

Here is the main theorem of this paper, which shows that the
density in the case of polygons is at most as in the case of
equilateral triangle.

\begin{theorem}
Under congruent Voronoi polygon conditions, the minimum density of
the two-coverage with equal circles is 2$\vartheta$.
\end{theorem}

\begin{figure}
\centering
\includegraphics[height=3.5cm]{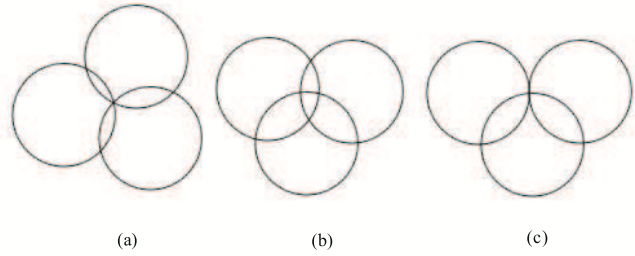}
\caption{(a) Only one point is shared by three circles. (b) An area
is shared by three circles. (c) The tangent point of two circles is
shared by three circles} \label{trangleMeet}
\end{figure}

\begin{figure}
\centering
\includegraphics[height=4.0cm]{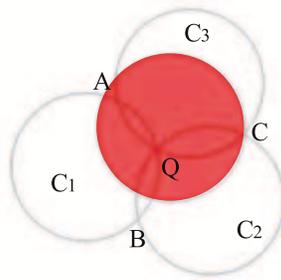}
\caption{A possible non-triangle tessellation derived by pattern
($a$) in Figure \ref{trangleMeet}.} \label{triangleNot}
\end{figure}

\begin{figure}
\centering
\includegraphics[height=4.8cm]{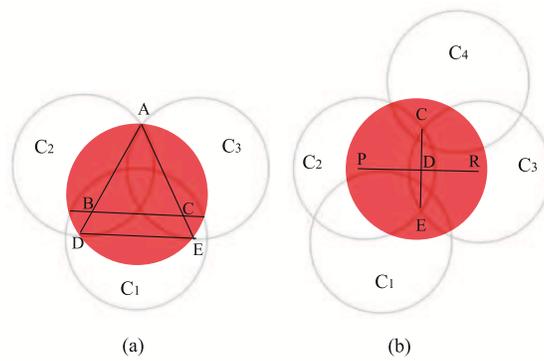}
\caption{Two types of two-coverage derived  by pattern ($b$) shown
in Figure \ref{trangleMeet}.} \label{triangleB}
\end{figure}

\begin{figure}
\centering
\includegraphics[height=6.7cm]{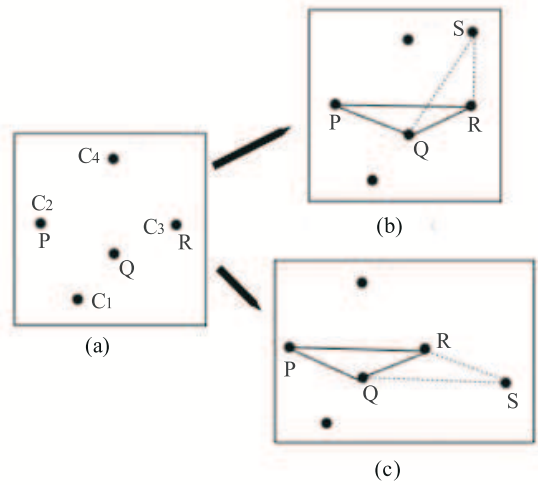}
\caption{ (a) The centers of five circles shown in Figure
\ref{triangleB}(b). (b) Point $S$ generated by rotating $\triangle
PQR$ (rotation symmetry tessellation). (c) Point $S$ generated by
mirror-mapping $\triangle PQR$ (mirror symmetry tessellation). }
\label{triangleBB}
\end{figure}

\begin{figure}
\centering
\includegraphics[height=4.3cm]{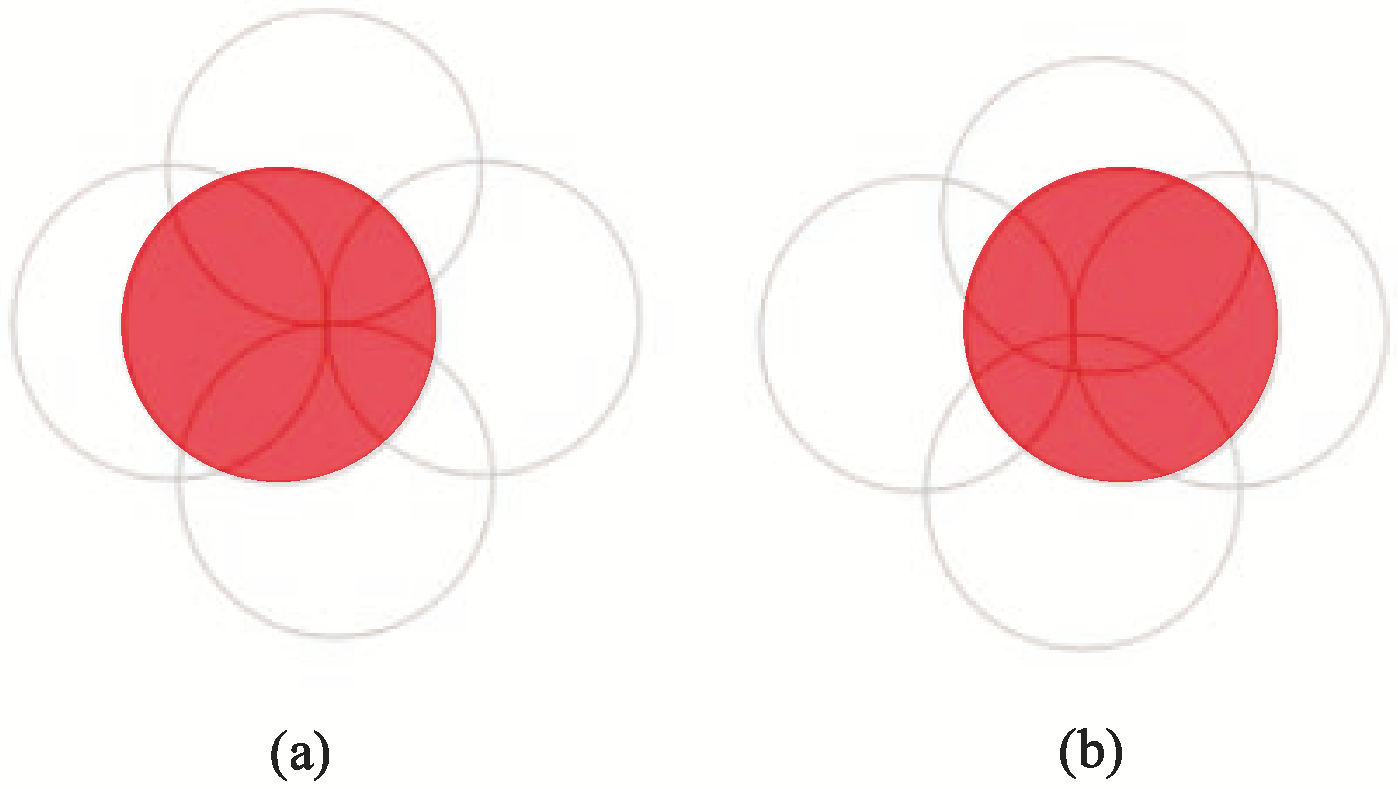}
\caption{ Two types of two-coverage derived by pattern ($c$) shown
in Figure \ref{trangleMeet}. (a) Two pairs of two circles are
tangent. (b) Only two circles are tangent.} \label{triangleC}
\end{figure}

\begin{proof}

Here we may assume that there are no coincident circles forming a
two-fold covering, since coincident cases cannot generate congruent
Voronoi polygons.

   To achieve two-coverage, there is at least one interior point of any circle,
   which belongs to at least three other circles, since any two circles
   cannot cover fully a circle unless that circle and one of them are coincident.
   A point that is shared by three circles may be classified into the following three cases:

\begin{description}
  \item[1)]As shown in Figure \ref{trangleMeet}(a), three circles meet at a point;
  \item[2)]As shown in Figure \ref{trangleMeet}(b), an area is shared by three circles;
  \item[3)]As shown in Figure \ref{trangleMeet}(c), two circles are tangent. Another circle passes through that tangent point.
\end{description}

\noindent Any tessellation pattern is one of the above three cases,
or their combination.

First we consider pattern ($a$) in Figure \ref{trangleMeet}. To
achieve two-coverage, it may result in two cases: (1) A circle is
fully covered by three other circles, i.e., a triangle pattern
shown in Figure \ref{triangleA}; (2) A circle is not fully covered
by three other circles, i.e., a non-triangle pattern shown in
Figure \ref{triangleNot}. Lemma 1 indicates that there exists indeed
such a triangle pattern. By Lemma 2, we obtain that the density of
the two-coverage in the case of the triangle pattern is at least
2$\vartheta$. For the non-triangle pattern, there must be an area
that shared by three circles, e.g., circles $C_1$, $C_2$ and red
circle shown in Figure \ref{triangleNot} constructs such an area.
Therefore, the non-triangle pattern is actually pattern ($b$) in
Figure \ref{trangleMeet}.

Now consider pattern ($b$) shown in Figure \ref{trangleMeet}. This
may derive two types of two-coverage: (1) As shown in Figure
\ref{triangleB}(a), the red circle is fully covered by circles
$C_1$, $C_2$, $C_3$; (2) As shown in Figure \ref{triangleB}(b), the
red circle requires at least four circles to be fully covered, i.e.,
circles $C_1$, $C_2$, $C_3$ cannot cover fully it. In Figure
\ref{triangleB}(a), we assume that the red circle passes through the
intersection $A$ of circles $C_2$, $C_3$. If it does not pass
through $A$, the density computation is similar. Under this
assumption, the Voronoi polygon lying in the red circle is a
triangle $ABC$, as depicted Figure \ref{triangleB}(a). Clearly, the
area of $\triangle ABC $ is no greater than that of the inscribed
triangle, i.e., $\triangle ADE$. From the fact that the area of a
regular triangle is maximum among all inscribed triangles, we can
obtain that the density of the tessellation pattern shown in Figure
\ref{triangleB}(a) is at least 2$\vartheta$.

   Compared with Figure \ref{triangleB}(a), computing the density of Figure
   \ref{triangleB}(b) is more difficult. Given the locations of a few
   centers, in many cases, the locations of the remaining centers can be
   computed by using the symmetry property, since under congruent Voronoi polygon conditions, there must be
   the symmetry property between centers of circles, i.e., each center
   has the same adjacent structure. In Figure \ref{triangleB}(b), there are five circles, the centers of which are shown in Figure
   \ref{triangleBB}(a). Now we compute another center $S$ by using three
   out of the five centers, i.e., $P$, $Q$, $R$. As shown in Figure \ref{triangleBB},
   we restrict ourselves to the case where centers $P$, $Q$, $R$ are not
   collinear because the case where centers $P$, $Q$, $R$ are collinear is trivial.
   Since center $R$ must have same adjacent structure as center $Q$, there must exist a triangle congruent to $\triangle
   PQR$. We can construct such a congruent triangle by only two
   operations: rotation operation and  mirror-map operation. Figure \ref{triangleBB}(b) depicts $\triangle QRS$ (dotted line) congruent to $\triangle
   PQR$, which is generated by a rotation operation. Figure \ref{triangleBB}(c) depicts $\triangle QRS$ (dotted line) congruent to $\triangle
   PQR$, which is generated by a mirror-map operation.
   Notice, with respect to line segment $\overline{PR}$, we can carry out repeatedly at
most $\frac{360}{\angle STR}$
   rotation operations. However, suppose that the plane to be covered may be arbitrarily large or infinite, then the number
   of rotation or mirror-map operations with respect to line segment $\overline{PR}$  should be able to be arbitrarily large or infinite.
    Therefore, with respect to the line segment $\overline{PR}$, there must exist mirror-map operations. It is easy to see that once a
   mirror-map operation is carried out, there must be another
   mirror-map to construct parallelogram $QRST$ congruent to parallelogram
   $PQRS$ shown in Figure \ref{triangleBB}(c), where centers $P$, $R$ and $T$ are collinear. The rest may be deduced by
   analogy. It implies that centers lying on the straight line through $\overline{PR}$ are evenly spaced,
   and the distance between any two adjacent centers on this line is $|PR|$.
   Similarly, this claim for straight line through $\overline{QS}$ holds also. Now we compute the density.
   Let the radius be 1, $|PR|=2x$, $y$ and $z$ be the distance from the centers
   of circles $C_4$ and $C_1$ to straight line $PR$, respectively.
   We may assume without loss of generality that $y\geq z$. Clearly,
   from Figure \ref{triangleB}(b), we have $y \leq |CD| + 1$, and
   $|CD|^2=1^2-|DR|^2 = 1 -x^2$. It implies $y \leq \sqrt{1-x^2} + 1$.
   The density of Figure \ref{triangleB}(b) is at least $\min(\frac{\pi}{xy},\frac{\pi}{xz}) = \frac{\pi}{xy} \geq \frac{\pi}{x(\sqrt{1-x^2} +
   1)}$. Notice,  centers based on line $PR$ yields at least the density of
   $\frac{\pi}{2xy}$. the density of centers based on line $QS$ is
   at least $\frac{\pi}{2xy}$  also.  When $x=\frac{\sqrt{3}}{2}$, the formula $x(\sqrt{1-x^2}+1)$
   reaches maximum, which is equal to $\frac{3\sqrt{3}}{4}$.
   Therefore, the density in this case is at least $\frac{4\pi}{3\sqrt{3}} =
   2\vartheta$. (In fact, we found that this proof approach is also
   suitable for the proof of Lemma 2, i.e., constructing congruent Voronoi
   triangles are based two straight lines, and centers on each straight line are evenly spaced.)

   Finally we consider pattern ($c$) shown in Figure \ref{trangleMeet}.
This pattern may derive two types of two-coverage shown in Figure
\ref{triangleC}. In Figure \ref{triangleC}(a), the red circle is
covered by four circles, the relation of which is that two pairs of
tangent circles meets at a point. In Figure \ref{triangleC}(b),
there is only one pair of tangent circles. Both the cases are an
extreme case of Figure \ref{triangleB}(b) with $xy \leq 1$. In a way
similar to Figure \ref{triangleB}(b) mentioned above, we can obtain
that in either case the density of pattern
   ($c$) is at least $\pi > 2\vartheta$, since the average distance between two
   adjacent centers in the horizontal direction is $r$, and the average distance between two
   adjacent centers in the vertical direction is $r$ in the case of
   Figure \ref{triangleC}(a), and at most $r$ in the case of Figure \ref{triangleC}(b), where $r$ is the radius of the circle.


 Gathering together the above discussion about three cases shown in Figure \ref{trangleMeet} completes the proof of the theorem.

\end{proof}

\tolerance=1
\emergencystretch=\maxdimen
\hyphenpenalty=10000
\hbadness=10000

In fact, any $k$-fold covering under congruent Voronoi polygon
conditions is a double lattice-covering, each of which has the same
structure, while Ref. \cite{Blundon:1957} is a single
lattice-covering. Therefore, we can prove Theorem 1 also by
extending some theorems of Ref.~\cite{Blundon:1957}.

\bibliographystyle{splncs}
\bibliography{two_cover}

\end{document}